\documentclass[11pt]{amsart}

\usepackage{amsmath,amsfonts,amsthm,amsopn,amssymb}
\usepackage{verbatim,wasysym,mathrsfs}
\usepackage[left=2.6cm,right=2.6cm,top=3.1cm,bottom=3.1cm]{geometry}
\usepackage{cite,marginnote}
\usepackage{color,enumitem,graphicx}
\usepackage[colorlinks=true,urlcolor=blue, citecolor=red,linkcolor=blue,
linktocpage,pdfpagelabels, bookmarksnumbered,bookmarksopen]{hyperref}
\usepackage[english]{babel}
\usepackage[T1]{fontenc}
\usepackage{lmodern}
\usepackage[hyperpageref]{backref}

\numberwithin{equation}{section}

\makeindex
\newtheorem{theorem}{Theorem}[section]
\newtheorem{proposition}[theorem]{Proposition}
\newtheorem{lemma}[theorem]{Lemma}
\newtheorem{remark}[theorem]{Remark}
\newtheorem{example}[theorem]{Example}
\newtheorem{corollary}[theorem]{Corollary}
\newtheorem{definition}[theorem]{Definition}

\newcommand{\s}{\section}

\newcommand{\R}{\mathbb R}

\newcommand{\bt}{\begin{theorem}}
\newcommand{\et}{\end{theorem}}
\newcommand{\bl}{\begin{lemma}}
\newcommand{\el}{\end{lemma}}
\newcommand{\bd}{\begin{definition}}
\newcommand{\ed}{\end{definition}}
\newcommand{\bc}{\begin{corollary}}
\newcommand{\ec}{\end{corollary}}
\newcommand{\bp}{\begin{proof}}
\newcommand{\ep}{\end{proof}}
\newcommand{\bx}{\begin{example}}
\newcommand{\ex}{\end{example}}
\newcommand{\bi}{\begin{exercise}}
\newcommand{\ei}{\end{exercise}}
\newcommand{\bo}{\begin{proposition}}
\newcommand{\eo}{\end{proposition}}
\newcommand{\br}{\begin{remark}}
\newcommand{\er}{\end{remark}}
\newcommand{\be}{\begin{equation}}
\newcommand{\ee}{\end{equation}}
\newcommand{\ba}{\begin{align}}
\newcommand{\ea}{\end{align}}
\newcommand{\bn}{\begin{enumerate}}
\newcommand{\en}{\end{enumerate}}
\newcommand{\bg}{\begin{align*}}
\newcommand{\bcs}{\begin{cases}}
\newcommand{\ecs}{\end{cases}}

\def\R{\mathbb R}

\def\N{\mathbb N}

\def\Proof{\noindent{\bf Proof}\quad}
\def\qed{\hfill$\square$\smallskip}
\makeatletter
\def\@makefnmark{}
\makeatother

\newcommand{\bean}{\begin{eqnarray*}}
\newcommand{\eean}{\end{eqnarray*}}

\renewcommand{\triangle}{\Delta}
\renewcommand{\epsilon}{\varepsilon}

\title[ Existence and multiplicity of bound state solutions]
{ Existence and multiplicity of bound state solutions to a Kirchhoff type equation
 with a general nonlinearity}

\author[Z. S. Liu]{Zhisu Liu}
\author[H. J. Luo]{Haijun Luo}
\author[J. J. Zhang]{Jianjun Zhang}

\address[Z. S. Liu]{\newline\indent Center for Mathematical Sciences,
\newline\indent
China University of Geosciences,
\newline\indent
Wuhan, Hubei, 430074, PR China}
\email{\href{mailto:liuzhisu183@sina.com}{liuzhisu183@sina.com}}

\address[H. J. Luo]{\newline\indent School of Mathematics,
\newline\indent
Hunan University,
\newline\indent
Changsha 410082, Hunan, P.R. China}
\email{\href{mailto:luohj@hnu.edu.cn}{luohj@hnu.edu.cn}}

\address[J. J. Zhang]{\newline\indent
College of Mathematics and Statistics,
\newline\indent
Chongqing Jiaotong University,
\newline\indent
 Chongqing 400074, PR China }
\email{\href{mailto:zhangjianjun09@tsinghua.org.cn}{zhangjianjun09@tsinghua.org.cn}}

\thanks{Z. Liu is supported by NSFC (No.11626127) and Hunan Natural Science Excellent Youth Fund (No.2020JJ3029). H. Luo is supported by the Fundamental Research Funds for the Central Universities(No.531118010205) and NSFC(No.11901182). J.\ Zhang was supported by NSFC(No.11871123).}


\begin{document}

\begin{abstract}
In this paper, we consider the following Kirchhoff type equation
$$
-\left(a+ b\int_{\R^3}|\nabla u|^2\right)\triangle {u}+V(x)u=f(u),\,\,x\in\R^3,
$$
where $a,b>0$ and $f\in C(\R,\R)$, and the potential $V\in C^1(\R^3,\R)$ is positive, bounded and satisfies suitable
decay assumptions. By using a new perturbation approach together with a new version of global compactness lemma
of Kirchhoff type, we prove the existence
and multiplicity of bound state solutions for the above problem with a general nonlinearity.
We especially point out that neither the corresponding Ambrosetti-Rabinowitz condition
nor any monotonicity assumption is required for $f$. Moreover, the potential $V$ may not be radially symmetry or coercive. As a prototype, the nonlinear term involves the power-type nonlinearity $f(u) = |u|^{p-2}u$ for $p\in (2, 6)$. In particular, our results generalize and improve the results by
Li and Ye (J.Differential Equations, 257(2014): 566-600), in the sense that the case $p\in(2,3]$ is left open there.
\end{abstract}
\maketitle
\begin{center}
\begin{minipage}{12cm}
\tableofcontents
\end{minipage}
\end{center}

\section{Introduction}
\label{sec1}

 In the present paper, we investigate the existence and multiplicity of bound state solutions to the following Kirchhoff equation
\begin{align}\label{K}\tag{K}
-\left(a+ b\int_{\R^3}|\nabla u|^2\right)\triangle {u}+V(x)u=f(u),\,\,x\in\R^3,\,\,u\in H^1(\R^3),
\end{align}
where $V\in C(\R^3,\R)$ and $a,b>0$ are positive constants.
Problem (K) arises in an interesting physical context.
Precisely, if we set $V(x)=0$ and a domain $\Omega\subset\R^3$ and replace $f(u)$ by $f(x,u)$,
problem (K) becomes as
the following Dirichlet problem:
\begin{equation} \label{eqn:youjielv}
\left\{
  \begin{aligned}
    -\left(a+b\int_{\Omega}|\nabla u|^2\right)\triangle {u} =f(x,u),\quad & \mbox{in}\,\,\Omega, \\
    u=0,\quad & \mbox{on}\,\, \partial\Omega,
   \end{aligned}
\right.
\end{equation}
which is the general form of the stationary counterpart of the hyperbolic Kirchhoff equation
\begin{equation}\label{eqn:stationary}
\rho\frac{\partial^2u}{\partial t^2}=\left[\frac{p_0}{h}
+\frac{E}{2L}\int_0^L\left(\frac{\partial u}{\partial x}\right)^2dx\right]\frac{\partial^2u}{\partial x^2}+f(t,x,u).
\end{equation}
This equation was proposed by Kirchhoff in \cite{Kirchhoff83} as an
existence of the classical D'Alembert's wave equations for free
vibration of elastic strings, and takes into account
the changes in length of the string produced by transverse
vibrations. In (\ref{eqn:stationary}), $L$ denotes the length of the string,
$E$ the Young modulus of the material, $h$ is the area of the cross
section, $\rho$ stands for mass density and $p_0$ is the initial tension, $f(t,x,u)$ stands for the external force.
The function $u$ denotes the displacement, $b$ is the initial tension while $a$ is related to the intrinsic properties.
Besides,
we also point out that Kirchhoff problems appear in other fields
like biological systems, such as population density, where $u$
describes a process which depends on the average of itself.
For the further physical background, we refer the readers to \cite{Cavalcanti01,Arosio96,Chipot97}.

\subsection{Overview and motivation}
Due to the presence of the integral term, Kirchhoff equations are no longer
a pointwise identity and therefore, are viewed as being nonlocal.
This fact brings mathematical challenges to the analysis, and meanwhile,
makes the study of such a problem particularly interesting. In the past decades, Kirchhoff problems have been receiving extensive attention. In particular, initiated by Lions \cite{Lions78},
the solvability of Kirchhoff
type equation (\ref{eqn:youjielv}) has been investigated in many studies,
 see \cite{Alves10,Alves05,Liang13,Ma03,Mao09,Perera06,Sun16,Zhang06,Shuai15,Tang16} and the references therein.

There also have been many interesting works about the existence and
multiplicity of bound state solutions to Kirchhoff type equation (K)
via variational methods, see for instance
\cite{Alves12,Azzollini11,He11,He11Zou,Li14,Figueiredo14,Liu15,Nie12,Wang12,Wu11,He16,Xie17,Liu17,Deng14,Xie19,Li20} and the references therein.
We note that minimax methods are used to study the existence and
multiplicity as a typical
way. In this process, one has to overcome the difficulties arising from
 the effect of non-local property and showing the boundedness and compactness of Palais-Smale ((PS) for short) sequences.
 For this aim, one usually assumes that the function $f$ satisfies either the $4$-superlinear
 growth condition:
\begin{align}\label{11111}\tag{4-superlinear}
\lim_{|t|\rightarrow+\infty}F(u)/u^4=+\infty,
\end{align}
where $F(u)=\int_{0}^uf(s)ds$, or the well-known Ambrosetti-Rabinowitz ((AR) for short) type condition
$$
0<F(u)\leq\frac{1}{\mu}f(u)u,\quad \mu>4,
$$
or the monotonicity condition
$$
\frac{f(u)}{u^3} \quad\text{is\, strictly\, increasing\, in\,\,} (0, +\infty).
$$
The above conditions are crucial in proving the existence and boundedness of (PS) sequences. Furthermore,
 nontrivial solutions can be obtained by providing some further conditions on $f$ and $V$ to guarantee
the compactness of the (PS) sequence, such as the radial symmetric setting or coercive condition.
It worth of pointing out that, without above conditions, Li and Ye \cite{Li14} proved the existence of positive ground state
solutions to problem (K) with $f(u)=|u|^{p-2}u$, $p\in(3,6)$ by using
the method of Nehari-Pohozaev manifold together with the concentration compactness arguments.
Recently, there results of \cite{Li14} were extended in \cite{Liu15} to the more general case, see also \cite{Guo15,Tang17}.

Compared with the existence results on nontrivial solutions,
there is few works published on the infinitely many solutions of
Kirchhoff type problem in $\R^3$, see \cite{Jin10,Wu11,Chen17,Nie12}.
As mentioned above, (AR)-condition or $4$-superlinear
 growth condition and some compactness conditions play important roles in this literatures.
 More specifically, Sun et al \cite{Sun19} obtained
infinitely many sign-changing solutions to problem (K) without $4$-superlinear
 growth condition but the coercive condition of $V$, by using a combination
of invariant sets and the Ljusternik-Schnirelman type minimax method.
Under
some weak compactness assumptions on $V$ without radial symmetry setting or compactness
hypotheses, Zhang et al. \cite{Zhang2020}
established the existence of infinitely many solutions to problem (K) with $f$ satisfying $4$-superlinear
 growth condition. Very recently, Liu et al. \cite{Liu19--} employed a novel perturbation approach and the method of invariant
sets of descending flow to
prove the existence of infinitely many sign-changing solutions to problem (K) with a general nonlinearity
in the radial symmetry setting.

\subsection{Our problem}
These results above left one question:
\begin{center}
{\it Does problem (K) admit infinitely many nontrivial solutions without the radial symmetric condition
or coercive condition in the case
$$f(u)\sim|u|^{p-2}u,\, p\in(2,4) ?$$}
\end{center}
Obviously, this type of nonlinearity $f$ does not satisfy (AR)-condition (or the $4$-superlinear)
or monotonicity assumptions mentioned as before.
To the best of our knowledge, so far there has been no results known in this aspect.
The main interest of the present paper is to give an affirmative answer to this question.
\vskip0.1in
\subsection{Our results} Throughout this paper, we assume nonlinearity $f$ satisfies the following hypotheses
\begin{itemize}
\item[ ($f_1$)]$f\in C(\R,\R)$ and $\lim\limits_{u\rightarrow0}\frac{f(u)}{u}=0$;
\item[ ($f_2$)]$\limsup\limits_{|u|\rightarrow\infty}\frac{|f(u)|}{|u|^{p-1}}<\infty$ for some $p\in(2,6)$;
\item[ ($f_3$)]there exists $\mu>2$ such that $uf(u)\geq\mu F(u)>0$ for $u\not=0$, where $F(u)=\int_{0}^uf(s)ds$.
\end{itemize}
These are quite natural assumptions when dealing with general subcritical nonlinearities.
In particular by ($f_1$)-($f_2$) it follows that
 for any $\epsilon>0$,
there exists $C_{\epsilon} > 0$ such that
\begin{equation}\label{eqn:fF}
|f(u)|\leq \epsilon |u|+C_{\epsilon}|u|^{p} \quad\text{and}\quad |F(u)| \leq \varepsilon u^{2}+C_{\varepsilon}|u|^{p+1}.
\end{equation}
\begin{remark}
It follows from $(f_1)$-$(f_3)$ that $2<\mu\leq p<6$. As a reference model,
$f(u) = |u|^{p-2}u$ satisfies ($f_1$)-($f_3$) for $p\in(2,6)$.
\end{remark}

Moreover, the potential $V\in C^1(\R^3,\R)$ enjoys
the following condition:
\begin{itemize}
\item[ ($V_1$)] there exist $V_0,V_1>0$ such that
$V_0\leq V(x)\leq V_1$ for all $x\in\R^3$;
\item[ ($V_2$)]for all $\gamma>0$, $\lim_{|x|\rightarrow\infty}
\frac{\partial V}{\partial r}(x)e^{\gamma|x|}=+\infty$, where $\frac{\partial V}{\partial r}(x)
=(\frac{x}{|x|},\nabla V(x))$;
\item[ ($V_3$)]there exists $\bar{c}>1$ such that $|\nabla V(x)|\leq \bar{c}\frac{\partial V}{\partial r}(x)$
for all $x\in\R^3$ and $|x|\geq \bar{c}$;
\item[ ($V_4$)] for all almost $x\in\R^{3}$, $(\nabla V(x), x)\in L^{\infty}(\R^{3})\cup L^{2}(\R^{3})$
and $\frac{\mu-2}{\mu}V(x)
\geq (\nabla V(x), x)\geq 0$.
\end{itemize}
\begin{remark}
We note that ($V_2$) and ($V_3$) were firstly given in Cerami et al \cite{Cerami05}
to study the existence of infinitely many bound state solutions for nonlinear scalar field equations.
This assumptions are key in recovering the compactness of solution sequence
when one uses local Pohozaev indentity together with decay estimates to study the behavior of solution,
see also Liu and Wang \cite{Liuj14}. Of course, ($V_4$) is also a very natural condition to ensure the boundedness of solution sequence, see Li and Ye \cite{Li14}.
It is not difficult to find some concrete function $V$ satisfying assumptions ($V_1$)-($V_4$), such as
$$
V(x)=V_1-\frac{1}{1+|x|},\quad V_1>\frac{(3\mu-2)C}{\mu-2}, \quad V_0\in(0,V_1-1)
$$
or
$$
V(x)=V_0+Ce^{\frac{-1}{1+|x|}},\quad V_0>\frac{\mu+2}{2\mu}Ce^{-1}, \quad V_1\in(V_0+C,+\infty),
$$
where $C$ is a positive constant.
\end{remark}
Our main result is as follows:

\begin{theorem}\label{Thm:existence}
If ($V_1$)-($V_4$) and ($f_1$)-($f_3$) hold, then problem (K)
admits at least one least energy solution in $H^1(\R^3)$.
\end{theorem}

\begin{theorem}\label{Thm:duochjie}
If ($V_1$)-($V_4$) and ($f_1$)-($f_3$) hold,
then problem (K) has infinitely many bound state solutions in $H^1(\R^3)$ provided that $f(u)$ is odd in $u$.
\end{theorem}
Now we summarize two main difficulties in finding bound state solutions to problem (K) under the effect of nonlocal
term $\int_{\R^3} |\nabla u|^2$.
On one hand, when $p\in(2,4)$, both the so-called 4-(AR) condition and the monotonicity condition fail,
which make tough to
get the boundedness of (PS) sequences.
On the other hand, it is also hard to prove the convergence of (PS) sequences without radial symmetry setting or compactness
hypotheses for $V$.
It is mainly motivated by \cite{Cerami05,Li14,Liuj14,Liuz16}
that we make use of a new perturbation approach together with
symmetric mountain-pass theorem to study problem (K). More precisely,
in order to get boundedness and compactness of (PS) sequences,
we modify problem (K) by adding a conceive term and a nonlinear
term growing faster than 4, see the modified problem (K$_\lambda$),
and then the corresponding Pohozaev type identity enables us to get a bounded solution sequence independent of the parameter $\lambda$.
As a result, by passing to the limit, a convergence argument allows us to get
nontrivial solutions of the original problem (K). In this process,
we also need to establish a version of global decomposition of solution sequences (may be containing sign-changing solutions)
which seems new for Kirchhoff type equations. This decomposition is crucial in using the local Pohozaev identity and some decay estimates of solutions
to prove compactness
of the sequence of solutions. Moreover, we believe that this perturbation approach should be of independent interest in other problems.

\begin{remark}
The first result is not surprising. Indeed, we can see \cite{Liu15}
where they proved the existence of positive ground states to problem (K)
with a general nonlinearity, and even some more general assumptions for $f$ were used in \cite{Figueiredo14,Guo15,Tang17} to
study the existence of ground state solutions.
However, the methods used in this paper are different from ones in \cite{Figueiredo14,Li14,Liu15,Guo15,Tang17}.
The core of this paper is proving the existence of infinite many solutions which seems nontrivial.
 But it seems difficult to obtain infinitely many solutions by using those arguments in \cite{Liu15,Li14,Figueiredo14,Guo15,Tang17}.
\end{remark}

Hereafter, the letter $C$ will be repeatedly used to denote various
positive constants whose exact values are irrelevant. We omit the symbol $dx$ in the integrals when no confusion
can arise.
This paper is organized as follows. Firstly, some
notations are given in Section 2, and Section 3 is devoted to the
existence of positive ground state solution. Then in Section 4, we
investigate the existence of infinitely many bound state solutions.

\s{Preliminary results}\label{sec:prel}

To proceed, we first define the Hilbert space
\begin{align*}
 H=\left\{u\in H^1(\R^3):\int_{\R^3} V(x)u^2<\infty\right\}
\end{align*}
with the inner product
$$
\langle u,v\rangle=\int_{\R^3}a\nabla
u\nabla v+V(x)uv
$$
and the norm
$$
\|u\|:=\sqrt{\langle u,u\rangle}=\left(\int_{\R^3}a|\nabla
u|^2+V(x)u^2\right)^{\frac{1}{2}}.
$$
The associated energy functional $I:H\rightarrow\R$ is given by
\begin{equation*}\label{eqn:fanhan}
I(u)=\frac{1}{2}\|u\|^2
+\frac{b}{4}\left(\int_{\R^3}|\nabla u|^2\right)^2-\int_{\R^3}F(u).
\end{equation*}
It is a well-defined $C^1$ functional in $H$ and its derivative is given by
\begin{equation*}\label{eqn:Frechet}
 I'(u)v=\int_{\R^3}(a\nabla u\nabla v+V(x)uv)+b\int_{\R^3}|\nabla u|^2\int_{\R^3}\nabla u\nabla v -\int_{\R^3} f(u)v,\,\,v\in H.
\end{equation*}
We introduce the following coercive function which will be of use
\begin{equation}\label{eqn:coercive0}
 W(x):=1+|x|^{\alpha},\quad 0<\alpha<\frac{\mu-2}{\mu},\,\,x\in\R^3.
\end{equation}
Obviously,
\begin{equation}\label{eqn:coercive}
W(x)\geq
1>0,\,\,\lim\limits_{|x|\rightarrow\infty}W(x)\rightarrow\infty,
\end{equation}
 and
\begin{equation}\label{eqn:differe}
 \frac{\mu-2}{\mu}W(x)\geq(\nabla W(x), x)\geq 0\quad\,\, x\in\R^{3}.
\end{equation}
Let $E_{\lambda}:=\{u\in H:\int_{\R^{3}}\lambda W(x)u^{2}{\rm d}x<\infty\}$
equipped with the norm
$$\|u\|_{E_\lambda}=\bigg(\int_{\R^{3}}(|\nabla u|^{2}+V(x)u^2+\lambda W(x)u^{2})\bigg)^{\frac{1}{2}}.$$
Note that $E=E_1\subset E_{\lambda}\subseteq H$ for $\lambda\in (0,1]$.

\section{Existence}\label{sec:subcritical}
\setcounter{equation}{0}

\subsection{The perturbed problem} It is known that the boundedness of the Palais-Smale sequence is not easy to prove for the case $p\in(2,4)$.
To overcome this difficulty, we introduce a perturbation technique to problem (K). We now give more details to
describe such a technique.
Fix $\lambda\in(0,1]$ and $r\in(\max\{p,4\}, 6)$, we consider the modified problem
$$
\left\{
  \begin{aligned}
-\left(a+ b\int_{\R^3}|\nabla u|^2\right)\triangle {u}+
V(x)u+\lambda W(x)u=f_{\lambda}(u),\quad\mbox{in}\,\,\R^3,\\
u\in E_\lambda,
  \end{aligned}
\right.\eqno{(K_{\lambda})}
$$
where
$$
f_{\lambda}(u)=f(u)+\lambda|u|^{r-2}u.
$$
An associated functional can be constructed as
$$
I_{\lambda}(u)=I(u)+\frac{\lambda}{2}\int_{\R^3}W(x)u^2-\frac{\lambda}{r}\int_{\R^3}|u|^r,\,\,u\in E_\lambda,
$$
and for $u,v\in E_\lambda$,
\begin{equation} \label{eqn:JJJ1}
I'_{\lambda}(u)v=\int_{\R^{3}}[a\nabla
u\nabla v+V(x)uv+\lambda W(x)uv]
+b\int_{\R^{3}}|\nabla u|^2\int_{\R^3}\nabla u\nabla v
-\int_{\R^{3}}(f(u)v+\lambda|u|^{r-2}uv).
\end{equation}
 It is known that $I_{\lambda}$ belongs to $C^{1}(E_\lambda,\R)$ or $C^{1}(E,\R)$ and a
critical point of $I_{\lambda}$ is a weak solution of problem
($K_\lambda$).
As we know, the original problem can be seen
as the limit system of ($K_\lambda$) as $\lambda\rightarrow 0^+$.

We will make use of the following Pohozaev type identity, whose proof is standard and
can be found in \cite{Berestycki1}.
\begin{lemma}\label{Lem:pohozaev}
Let $u$ be a critical point of $I_{\lambda}$ in $E_{\lambda}$ for $\lambda\in
(0,1]$, then
$$
\aligned
\frac{a}{2}\int_{\R^3}|\nabla u|^2+&\frac{3}{2}\int_{\R^N}(V(x)+\lambda W(x))u^2+\frac{1}{2}\int_{\R^3}(\nabla
V(x),x)u^2+\frac{\lambda}{2}\int_{\R^3}(\nabla
W(x),x)u^2\\+&\frac{b}{2}\left(\int_{\R^3}|\nabla u|^2\right)^2
-3\int_{\R^3}(F(u)+\frac{\lambda}{r}|u|^r)=0.
\endaligned
$$
\end{lemma}
We now verify that the functional $I_\lambda$ has the Mountain Pass geometry
uniformly in $\lambda$.

\begin{lemma}\label{Lem:MP1}
Suppose that ($V_1$)-($V_4$) hold.
 Then
 \begin{itemize}
 \item [(1)] there exist $\rho,\delta>0$ such that, for any $\lambda\in(0,1]$,
  $I_{\lambda}(u)\geq \delta$ for every $u\in S_\rho=\{u\in E_\lambda: \|u\|_{E_\lambda}=\rho\}$;
 \item [(2)] there is $v\in E\setminus\{0\}$ with $\|v\|_{E_\lambda}>\rho$ such that,
 for any $\lambda\in(0,1]$, $I_{\lambda}(v)<0$.
 \end{itemize}
\end{lemma}
\Proof
(1)  For any $u\in E_\lambda$, by the definition of $I_\lambda$, (\ref{eqn:fF})
and Sobolev's inequaity, one has
\begin{equation*}
\aligned
I_{\lambda}(u)
&\geq\frac{1}{4}\|u\|_{E_\lambda}^2-C\int_{\R^3}|u|^{p}-\frac{1}{r}\int_{\R^3}|u|^{r} \\
&\geq\frac{1}{4}\|u\|_{E_\lambda}^2
-C\|u\|_{E_\lambda}^{p}-\frac{C}{r}\|u\|_{E_\lambda}^{r}.
\endaligned
\end{equation*}
Taking $\rho>0$ small
enough, it is easy to check that there exists $\delta>0$ such that
$I_{\lambda}(u)\geq \delta$ for every $u\in S_{\rho}$. \\

(2) For $e\in E\setminus\{0\}$,
let $e_{t}=t^{1/2}e(\frac{x}{t})$. Observe that
\begin{equation*}\label{eqn:MPG2}
\int_{\R^3}F(e_{t})=t^{3}\int_{\R^3}F(t^\frac{1}{2} e)=:t^{3} \Phi(t).
\end{equation*}
By ($f_3$), a straightforward computation yields
$$
\frac{\Phi'(t)}{\Phi(t)}\geq \frac{\mu}{2t}, \quad \forall t>0
$$
and then, by integrating on $[1, t]$, with $t>1$, we have
$\Phi(t)\geq \Phi(1)t^{\frac{\mu}{2}}$, implying that
\begin{equation}\label{eq:F}
\int_{\mathbb R^{3}} F(e_{t}) \geq t^{\frac{\mu+6}{2}}\int_{\mathbb R^{3}}F(e).
\end{equation}
Then
by the definition of $I_\lambda$ and (V$_1$) and (\ref{eqn:coercive0}), one has
\begin{equation}\label{eqn:MPG}
\aligned
I_{\lambda}(e_{t})
&<\frac{t^2}{2}\|\nabla e\|_{2}^2
+\frac{t^4}{4}\|\nabla e\|_{2}^4+\frac{t^4}{2}\int_{\R^3}V(tx)e^2+\frac{\lambda t^4}{2}\int_{\R^3}W(tx)e^2
-t^{\frac{\mu+6}{2}}\int_{\R^3}F(e)\\
&\leq\frac{t^2}{2}\|\nabla e\|_{2}^2
+\frac{t^4}{4}\|\nabla e\|_{2}^4+\frac{t^4V_1}{2}\int_{\R^3}e^2+\frac{t^{4+\alpha}}{2}\int_{\R^3}W(x)e^2
-t^{\frac{\mu+6}{2}}\int_{\R^3}F(e)\\
&<0,
\endaligned
\end{equation}
which holds for $t>1$ large enough, owing to $\alpha<\frac{\mu-2}{\mu}$. The proof is complete.
\qed

By recalling the well-known Mountain-Pass theorem (see \cite{AR,Willem96}),
there exists a $(PS)_{c_{\lambda}}$ sequence
$\{u_n\}\subset E_\lambda$, that is,
\begin{equation*}
I_{\lambda}(u_n)\rightarrow c_{\lambda}\quad\text{and}\quad I'_{\lambda}(u_n)\rightarrow0.
\end{equation*}
We stress that $\{u_{n}\}$ depends on $\lambda$ but we omit this dependence in the sequel for convenience.
Here $c_{\lambda}$ is the Mountain Pass level characterized by
\begin{equation*}
c_{\lambda}=\inf\limits_{\gamma\in\Gamma_{\lambda}}\max\limits_{t\in[0,1]}I_{\lambda}(\gamma(t))
\end{equation*}
with
$$
\Gamma_{\lambda}:=\left\{\gamma\in C^1([0,1],E_\lambda):\,\gamma(0)=0\quad\text{and}\quad I_{\lambda}(\gamma(1)) <0\right\}.
$$
\begin{remark}\label{rem:bounded}
Observe from Lemma \ref{Lem:MP1} that there exist two constants $m_{1}, m_{2}>0 $
independently on $\lambda$ such that $m_{1}<c_\lambda<m_{2}$.
\end{remark}

In what follows, we prove the functional $I_{\lambda}$ satisfies the (PS)-condition.

\begin{lemma}\label{Lem:PS}
Assume that there exists $\{u_n\}\subset E_{\lambda}$ such that $I_{\lambda}(u_n)\rightarrow c_\lambda$ and
$I'_{\lambda}(u_n)\rightarrow 0$ for any fixed $\lambda\in(0,1)$ as $n\rightarrow\infty$,
then there exists a convergence subsequence of $\{u_n\}$,
still denoted by $\{u_n\}$, such that $u_n\rightarrow u$ in $E_\lambda$ for some $u\in E_\lambda$.
\end{lemma}
\Proof
For $\gamma\in(4,r)$, by (\ref{eqn:fF}) we have
$$
\aligned
&\gamma I_{\lambda}(u_n)-\langle I'_{\lambda}(u_n),u_n\rangle\\
&=\frac{\gamma-2}{4}\|u_n\|_{E_\lambda}^2+\frac{b(\gamma-4)}{4}\left(\int_{\R^3}|\nabla u_n|^2\right)^2\\
&+\int_{\R^3}\bigg(\gamma f(u_n)u_n-F(u_n)\bigg)+\lambda\frac{r-\gamma}{r}\int_{\R^3}|u|^r.
\endaligned
$$
Then it follows from (\ref{eqn:fF}) that
\begin{equation}\label{eqn:T5-}
\|u_n\|_{E_\lambda}^2+b\left(\int_{\R^3}|\nabla u_n|^2\right)^2+\lambda\int_{\R^3}|u|^r\leq C(1+\|u_n\|_{E_\lambda}+\|u_n\|_p^p)
\end{equation}
for large $n$. We claim that $\{u_n\}$ is uniformly bounded in $E_\lambda$. Assume by contradiction that $\|u_n\|_{E_\lambda}\rightarrow\infty$, then by (\ref{eqn:T5-}) we have
\begin{equation}\label{eqn:T6-}
\|u_n\|_{E_\lambda}^2+b\left(\int_{\R^3}|\nabla u_n|^2\right)^2+\lambda\|u_n\|_r^r\leq C\|u_n\|_p^p,
\end{equation}
which implies that
$$
\|u_n\|_2^2+\|u_n\|_r^r\leq C\|u_n\|_p^p.
$$
Let $t\in(0,1)$ be such that $\frac{1}{p}=\frac{t}{2}+\frac{1-t}{r}$. From the interpolation inequality,
we deduce that
\begin{equation}\label{eqn:T7-}
\|u\|_2^2+\|u\|_r^r\leq C \|u_n\|_p^p\leq C \|u_n\|_2^{p t}\|u_n\|_r^{p(1-t)}.
\end{equation}
It follows from (\ref{eqn:T7-}) that there exist $C_{1}$, $C_{2}>0$ such that
\begin{equation}\label{eqn:T8-}
C_{1}\|u_n\|_2^{\frac{2}{r}}\leq \|u_n\|_r\leq C_{2}\|u_n\|_2^{\frac{2}{r}}.
\end{equation}
In view of (\ref{eqn:T7-}) and (\ref{eqn:T8-}), we have $\|u_n\|_p^p\leq C_{3}\|u_n\|_2^2$ for some
$C_{3}>0$. Therefore, by (\ref{eqn:T6-}), we have for some $C_{4}>0$ such that
\begin{equation*}\label{eqn:T9-}
\|u_n\|_{E_\lambda}^2+b\left(\int_{\R^3}|\nabla u_n|^2\right)^2+\lambda\|u_n\|_r^r\leq C_{4}\|u_n\|_2^2.
\end{equation*}
Let $v_n=\frac{u_n}{\|u_n\|_{E_\lambda}}$, then
\begin{equation}\label{eqn:T9--}
\|v_n\|_2^2\geq \frac{1}{C_{4}}
\end{equation}
and
\begin{equation*}\label{eqn:T10-}
b\left(\int_{\R^3}|\nabla v_n|^2\right)^2\leq C_{4}\|u_n\|_{E_\lambda}^{-2},
\end{equation*}
which implies that $\int_{\R^3}|\nabla v_n|^2\rightarrow0$ as $n\rightarrow\infty$. By $\|v_n\|_{E_\lambda}=1$, we assume
$v_n\rightharpoonup v$ in $E_\lambda$. By Fatou's lemma we have
$$
\int_{\R^3}|\nabla v|^2\leq \liminf\limits_{n\rightarrow\infty}\int_{\R^3}|\nabla v_n|^2=0,
$$
which implies $v=0$. Then by (\ref{eqn:T9--}) we have $\|v\|_2^2\geq \frac{1}{C_{4}}$,  a contradiction.
Thus, we finish the proof of the claim.
 Without loss of generality, we assume that there exists $u\in E_\lambda$ such that
$$
\aligned
&u_n\rightharpoonup u\,\,\text{weakly\,in}\,\,E_\lambda,\\
& u_n\rightarrow u\,\,\text{strongly\,in}\,\,L^q(\R^3)\,\,\text{for}\,\,q\in[2,6).
\endaligned
$$
Note that
\begin{equation}\label{eqn:PS1}
\aligned
&(I'_{\lambda}(u_n)-I'_{\lambda}(u))(u_n-u)\\
&=\|u_n-u\|_{E_\lambda}^2+b\int_{\R^3}|\nabla u_n|^2\int_{\R^3}|\nabla(u_n-u)|^2\\
& +b(\int_{\R^3}|\nabla u_n|^2-\int_{\R^3}|\nabla u|^2)\int_{\R^3}\nabla u\nabla(u_n-u)
-\int_{\R^3}(f(u_n)-f(u))(u_n-u)\\
&-\lambda\int_{\R^3}(|u_n|^{r-2}u_n-|u|^{r-2}u)(u_n-u).
\endaligned
\end{equation}
According to the boundedness of $\{u_n\}$ in $E_\lambda$, one has
$$
\aligned
&b(\int_{\R^3}|\nabla u_n|^2-\int_{\R^3}|\nabla u|^2)\int_{\R^3}\nabla u\nabla(u_n-u)\rightarrow0.
\endaligned
$$
 Similarly, we also have
$$
\aligned
&\int_{\R^3}(f(u_n)-f(u))(u_n-u)\rightarrow0,\\
&\lambda\int_{\R^3}(|u_n|^{r-2}u_n-|u|^{r-2}u)(u_n-u)\rightarrow0,\quad \text{as}\,\,n\rightarrow\infty.
\endaligned
$$
Based on the above facts, from (\ref{eqn:PS1}) we deduce that $u_n\rightarrow u$ in $E_\lambda$.
\qed

It follows from Lemma \ref{Lem:PS}
that for each $\lambda\in(0,1]$, there exists $u_\lambda\in E_\lambda$ such that
$$
I_\lambda(u_\lambda)=c_\lambda \quad \text{and}\quad I'_\lambda(u_\lambda)=0.
$$
That is to say, $u_\lambda$ is a nontrivial solution of ($K_\lambda$).
 We now expect that $\{u_\lambda\}$
converges to a nontrivial solution of (K) as $\lambda\rightarrow0$ by controlling $\{u_\lambda\}$ in a proper way.

\begin{lemma}\label{Lem:local}
Suppose that $\lambda_n\rightarrow0^+$ as $n\rightarrow\infty$, $\{u_n\}\subset E_{\lambda_n}$ are nontrivial solutions of
($K_{\lambda_n}$) with $|I_{\lambda_n}(u_n)|\leq C$. Then there exists $M>0$ such that
$\|u_n\|_{E_{\lambda_n}}\leq M$ for some $M>0$ independently of $n$, and, up to subsequence, there is a solution $u\in H$ such that $u_n\rightharpoonup u_0$ in $H$.
\end{lemma}
\Proof
By sequence $\{\lambda_n\}\subset(0,1]$ satisfying $\lambda_n\rightarrow0^+$,
we can find a subsequence of $\{u_{\lambda_n}\}$ (still denoted by $\{u_n\}$)
of $I_{\lambda_n}$ with
$I_{\lambda_n}(u_n)=c_{\lambda_n}$. We claim that $\{u_n\}$ is bounded in $H$.
By the conditions of this lemma, we have
\begin{equation}\label{eqn:fun1}
\aligned
C\geq I_{\lambda_n}(u_n)=\frac{a}{2}\int_{\R^3}|\nabla u_n|^2&+\frac{1}{2}\int_{\R^3}(V(x)+\lambda_n W(x))u_n^2
\\
&+\frac{b}{4}\left(\int_{\R^3}|\nabla u_n|^2\right)^2-\int_{\R^3}F(u_n)-\frac{\lambda_n}{r}\int_{\R^3}|u_n|^r
\endaligned
\end{equation}
and
\begin{equation}\label{eqn:fun2}
\aligned
0=a\int_{\R^3}|\nabla u_n|^2&+\int_{\R^3}(V(x)+\lambda_n W(x))u_n^2\\
&+b\left(\int_{\R^3}|\nabla u_n|^2\right)^2-\int_{\R^3}f(u_n)u_n-\lambda_n\int_{\R^3}|u_n|^r.
\endaligned
\end{equation}
Moreover, from Lemma \ref{Lem:pohozaev},  the following identity holds
\begin{equation}\label{eqn:fun3}
\aligned
& \frac{a}{2}\int_{\R^3}|\nabla u_n|^2+\frac{3}{2}\int_{\R^3}(V(x)+\lambda_nW(x))u_n^2+\frac{1}{2}\int_{\R^3}(\nabla
V(x)+\lambda_n\nabla W(x),x)u_n^2\\+&\frac{b}{2}\left(\int_{\R^3}|\nabla u_n|^2\right)^2
-3\int_{\R^3}(F(u_n)+\frac{\lambda_n}{r}|u_n|^r)=0.
\endaligned
\end{equation}
Multiplying (\ref{eqn:fun1}), (\ref{eqn:fun2}) and (\ref{eqn:fun3}) by $4$, $-\frac{1}{\mu}$ and  $-1$ respectively and adding them up,
we get
\begin{equation*}\label{eqn:fun4}
\aligned
4C\geq&a\frac{3\mu-2}{2\mu}\int_{\R^3}|\nabla u_n|^2
+\frac{\mu-2}{2\mu}\int_{\R^3}(V(x)+\lambda_n W(x))u_n^2
-\frac{1}{2}\int_{\R^3}(\nabla
V(x)+\lambda_n\nabla W(x),x)u_n^2\\
&+\frac{\mu-2}{2\mu}b\left(\int_{\R^3}|\nabla u_n|^2\right)^2
+\lambda_n\frac{r-\mu}{\mu r}\int_{\R^3}|u_n|^r+\int_{\R^3}(\frac{1}{\mu}f(u_n)u_n-F(u_n)).
\endaligned
\end{equation*}
It then follows from ($V_4$) and (\ref{eqn:differe}) that
\begin{equation*}\label{eqn:fun5}
4C\geq a\frac{3\mu-2}{2\mu}\int_{\R^3}|\nabla u_n|^2+\frac{\mu-2}{2\mu}b\left(\int_{\R^3}|\nabla u_n|^2\right)^2+\lambda_n\frac{r-\mu}{\mu r}\int_{\R^3}|u_n|^r,
\end{equation*}
which implies that there exists $C_5>0$ independent of $\lambda_n$ such that
\begin{equation}\label{eqn:fun6}
\int_{\R^3}|\nabla u_n|^2<C_5.
\end{equation}
Moreover, combining (\ref{eqn:fF}), (\ref{eqn:fun1}) and hypotheses (V$_1$), we obtain that
for small $\varepsilon>0$, there exists $C_\varepsilon>0$ such that
\begin{equation}\label{eqn:fun7}
\aligned
C&>\frac{a}{2}\int_{\R^3}|\nabla u_n|^2+\frac{1}{2}\int_{\R^3}(V(x)+\lambda_n W(x))u_n^2
-\int_{\R^3}F(u_n)-\frac{\lambda_n}{r}\int_{\R^3}|u_n|^r\\
&>\frac{1-\varepsilon}{2}\int_{\R^3}V(x)u_n^2-C_\varepsilon\int_{\R^3}u_n^6+\frac{\lambda_n}{2}\int_{\R^3} W(x)u_n^2\\
&>\frac{1-\varepsilon}{2}\int_{\R^3}V(x)u_n^2-C_\varepsilon S^{-3}\left(\int_{\R^3}|\nabla u_n|^2\right)^3+\frac{\lambda_n}{2}\int_{\R^3} W(x)u_n^2.
\endaligned
\end{equation}
Combining (\ref{eqn:fun6}) and (\ref{eqn:fun7}),
there exists $C_6>0$ independent of $\lambda_n$ such that
\begin{equation}\label{eqn:fun7+}
\int_{\R^3}|\nabla u_n|^2+\int_{\R^3}(V(x)+\lambda_n W(x))u_n^2\leq C_6.
\end{equation}
The conclusions follow immediately.
\qed

The following lemma is devoted to the behavior of solution sequence to problem ($K_\lambda$).
\begin{lemma}\label{Lem:deco}
Let $\{u_n\}\subset E_\lambda$ be a solution sequence of problem ($K_\lambda$)
with $\lambda=\lambda_n\geq0$ and $\lambda_n\rightarrow0$, and
$\|u_n\|_{E_{\lambda_n}}\leq M$ for $M>0$ independent of $n$. Then there exist a
subsequence of $\{u_n\}$, still denoted by $\{u_n\}$, a number
$k\in\N\cup\{0\}$, and finite sequences
 $$
 (a_1,...,a_k)\subset \R,\quad (u_0,w_1,...,w_k)\subset H,\quad a_j\geq0,\,\,w_j\not\equiv0,
 $$
and $A\geq0$ and
$k$ sequences of points $\{y_n^j\}\subset\R^3$, $1\leq j\leq k$,
such that
\begin{itemize}
     \item[\rm (i)] $u_n\rightharpoonup u_0$, $u_n(\cdot+x_n^j)\rightharpoonup w_j$ in $H$ as $n\rightarrow\infty$,
    \item[\rm (ii)] $|y_n^j|\rightarrow+\infty$, $|y_n^j-y_n^i|\rightarrow+\infty\quad$if $i\neq j,n\rightarrow+\infty$,
    \item[\rm (iii)]   $\|u_n-u_0-\sum_{i=1}^{k}w_{i}(\cdot-y_n^{i})\|\rightarrow0$,
    \item[\rm (iv)]  $A=\|\nabla u_0\|_2^2+\sum_{i=1}^{k}\|\nabla w_{i}\|_2^2$,
    \item[\rm (v)]   for any $\varphi\in C_0^\infty(\R^3)$ with $\varphi\geq0$
    \begin{equation}\label{eqn:deco0}
(a+bA)\int_{\R^3}\nabla|w_j|\nabla\varphi+(V_0+a_j)\int_{\R^3}|w_j|\varphi
\leq \int_{\R^3}|f(w_j)|\varphi.
\end{equation}
\end{itemize}
\end{lemma}
\Proof
Note that $\{u_n\}$ is a bounded sequence in $H$. There exists $u_0\in H$ and $A>0$
such that $u_n\rightharpoonup u_0$ weakly in $H$ and $\|\nabla u_n\|_2^2\rightarrow A$ as $n\rightarrow\infty$ after extracting a subsequence.
For
 any $\psi\in C_0^\infty(\R^3)$, we have $J'_{\lambda_n}(u_n)\psi\equiv0$, where
 $$
 J_\lambda(u):=\frac{1}{2}\|u\|^2+\frac{\lambda}{2}\int_{\R^3}W(x)u^2
+\frac{Ab}{2}\int_{\R^3}|\nabla u|^2-\int_{\R^3}F(u)-\frac{\lambda}{r}\int_{\R^3}|u|^r.
 $$
 Moreover, one has for any $\psi\in C_0^\infty(\R^3)$
$$
\bigg|\lambda_{n}\int_{\R^{3}}W(x)u_{n}\psi\bigg|\leq\bigg(\lambda_{n}\int_{\R^{3}}W(x)u_{n}^{2}\psi\bigg)^{\frac{1}{2}}
\bigg(\lambda_{n}\int_{\R^{3}}W(x)\psi\bigg)^{\frac{1}{2}}\leq C\lambda_{n}^{\frac{1}{2}}\rightarrow 0,
$$
which, together with the fact that $J'_{\lambda_n}(u_n)=0$, implies that
$$
\lim\limits_{n\rightarrow\infty}J'(u_n)\psi=\lim\limits_{n\rightarrow\infty}
\left(J'_{\lambda_n}(u_n)\psi-\lambda_{n}\int_{\R^{3}}W(x)u_{n}\psi
+\lambda_n\int_{\R^3}|u_n|^{r-2}u_n\psi\right)=0,
$$
where the functional $J=J_\lambda$ with $\lambda=0$.
It then follows that $J'(u_0)=0$, that is,
\begin{equation}\label{eqn:deco1-00}
\int_{\R^3}(a\nabla u_0\nabla\psi+V(x)u_0\psi)+bA\int_{\R^3}\nabla u_0\nabla\psi
= \int_{\R^3}f(u_0)\psi.
\end{equation}
We \textbf{\emph{claim}} that the following differential inequality holds for any $\varphi\in C_0^\infty(\R^3)$ with $\psi\geq0$
\begin{equation}\label{eqn:deco1-0+}
\int_{\R^3}(a\nabla|u_0|\nabla\varphi+V_0|u_0|\varphi)+bA\int_{\R^3}\nabla|u_0|\nabla\varphi
\leq \int_{\R^3}|f(u_0)|\varphi.
\end{equation}
Set $u_\epsilon=\sqrt{|u_0|^2+\epsilon^2}-\epsilon$, $\epsilon>0$. It is clear that $u_\epsilon\rightarrow |u_0|$ in
$H$ as $\epsilon\rightarrow0$.
By (\ref{eqn:deco1-00}), we have for $\varphi\in C_0^\infty(\R^3)$ with $\varphi\geq0$
\begin{equation}\label{eqn:deco-2}
\aligned
(a+bA)\int_{\R^3}\nabla u_\epsilon\nabla\varphi&=
(a+bA)\int_{\R^3}\frac{u_0}{(|u_0|^2+\epsilon^2)^{\frac{1}{2}}}\nabla u_0\nabla\varphi\\
&=(a+bA)\bigg(\int_{\R^3}\nabla u_0\nabla\bigg(\frac{u_0\varphi}{(|u_0|^2+\epsilon^2)^{\frac{1}{2}}}\bigg)
-\int_{\R^3}|\nabla u_0|^2\frac{\epsilon^2\varphi}{(|u_0|^2+\epsilon^2)^{\frac{3}{2}}}\bigg)\\
&\leq (a+bA)\int_{\R^3}\nabla u_0\nabla\bigg(\frac{u_0\varphi}{(|u_0|^2+\epsilon^2)^{\frac{1}{2}}}\bigg)\\
&=-\int_{\R^3}V(x)u_0\frac{u_0\varphi}{(|u_0|^2+\epsilon^2)^{\frac{1}{2}}}
+\int_{\R^3}\frac{ f(u_0)u_0\varphi}{(|u_0|^2+\epsilon^2)^{\frac{1}{2}}}.
\endaligned
\end{equation}
So, from ($V_1$) we deduce that
\begin{equation}\label{eqn:deco-3-0}
a\int_{\R^3}\nabla u_\epsilon\nabla\varphi\leq-\int_{\R^3}V_0\frac{|u_0|^2\varphi}{(|u_0|^2+\epsilon^2)^{\frac{1}{2}}}
-bA\int_{\R^3}\nabla u_\epsilon\nabla\varphi+\int_{\R^3}\frac{f(u_0)u_0\varphi}{(|u_0|^2+\epsilon^2)^{\frac{1}{2}}}.
\end{equation}
Let $\epsilon\rightarrow0$ in (\ref{eqn:deco-3-0}), we obtain (\ref{eqn:deco5})
for $\varphi\in H^1(\R^3)$ with $\varphi\geq0$. The \textbf{\emph{claim}} is true.
We now apply the concentration compactness principle to the sequence of
 $\{v_{1,n}\}$ with $v_{1,n}=u_n-u_0$. Clearly, $v_{1,n}\rightharpoonup0$ weakly in $H$. If vanishing occurs,
 $$
 \sup\limits_{y\in\R^3}\int_{B_1(y)}|u_n-u_0|^2dx\rightarrow0,\quad\text{as}\,\,n\rightarrow\infty.
 $$
Then $v_{1,n}\rightarrow0$ in $L^s(\R^3)$
for $s\in(2,6)$. By the fact that $J'(u_0)=J_{\lambda_n}'(u_n)=0$, we arrive at
$$
\aligned
(a+bA)\int_{\R^3}|\nabla u_0|^2+\int_{\R^3}V(x)u_0^2&\leq
\liminf\limits_{n\rightarrow\infty}\bigg((a+bA)\int_{\R^3}|\nabla u_n|^2+\int_{\R^3}V(x)u_n^2\bigg)\\
&\leq\limsup\limits_{n\rightarrow\infty}\bigg((a+bA)\int_{\R^3}|\nabla u_n|^2+\int_{\R^3}V(x)u_n^2
+\lambda_n\int_{\R^3}W(x)u_n^2\bigg)\\
&=\limsup\limits_{n\rightarrow\infty}\bigg(\int_{\R^3}f(u_n)u_n+\lambda_n\int_{\R^3}|u_n|^r\bigg)\\
&\leq\int_{\R^3}f(u_0)u_0=(a+bA)\int_{\R^3}|\nabla u_0|^2+\int_{\R^3}V(x)u_0^2,
\endaligned
$$
which implies that $u_n\rightarrow u_0$ strongly in $H$.
So the conclusions of Lemma \ref{Lem:deco} hold for $k=0$.
If non-vanishing occurs,
then there exist $m>0$ and a sequence $\{y^1_n\}\subset\R^3$ such that
\begin{equation}\label{eqn:4.04}
\liminf\limits_{n\rightarrow\infty}\int_{B_{1}(y_n^1)}|v_{1,n}(x)|^2\geq m>0.
\end{equation}
 Let us consider
the sequence $\{v_{1,n}(\cdot+y_n^1)\}$. The boundedness of $\{v_{1,n}\}$ in $H$ implies
that there exists $w_1$ such that $v_{1,n}(\cdot+y_n^1)\rightharpoonup w_1$
in $H$. Furthermore, by (\ref{eqn:4.04}) one has
$$
\int_{B_{1}(0)}|w_1(x)|^2>\frac{m}{2},
$$
and, thus, $w_1\neq0$. Recalling the fact that $v_{1,n}\rightharpoonup 0$ in $H$, we know that $\{y_n^1\}$
must be unbounded and, up to a subsequence, we suppose that $|y^1_n|\rightarrow+\infty$.

Now we show the following inequality holds:
\begin{equation}\label{eqn:deco5}
(a+bA)\int_{\R^3}\nabla|w_1|\nabla\psi+\int_{\R^3}(a_1+V_0)|w_1|\psi
\leq \int_{\R^3}|f(w_1)|\psi
\end{equation}
for $\psi\in C_0^\infty(\R^3)$ with $\psi\geq0$.
Recalling (\ref{eqn:fun7+}), we have $\lambda_n\int_{\R^3}W(x)u_n^2\leq C$. So, (\ref{eqn:4.04}) implies that
$$
\aligned
C&\geq\lambda_n\int_{\R^3}W(x)|v_{1,n}(x)|^2\\
&\geq \lambda_n W(y_n^1)\int_{B_1(y_n^1)}|v_{1,n}(x)|^2-\lambda_n\int_{B_1(y_n^1)}|W(x)-W(y_n^1)||v_{1,n}(x)|^2\\
&\geq \lambda_n W(y_n^1)m-\lambda_n C,
\endaligned
$$
which implies that, up to subsequence, $\lambda_n W(y_n^1)\rightarrow a_1\in[0,+\infty)$.
Based on the above facts, we have for $\psi\in C_0^\infty(\R^3)$ with $\psi\geq0$
\begin{equation}\label{eqn:deco6}
\aligned
&\lambda_n\int_{\R^3}W(x+y_n^1)v_{1,n}(x+y_n^1)\psi\\
&=\lambda_n W(y_n^1)\int_{\R^3}v_{1,n}(x+y_n^1)\psi
+\lambda_n\int_{\R^3}(W(x+y_n^1)-W(y_n^1))v_{1,n}(x+y_n^1)\psi\\
&=a_1\int_{\R^3}v_{1,n}(x+y_n^1)\psi+o(1)\\
&=a_1\int_{\R^3}w_{1}\psi+o(1).
\endaligned
\end{equation}
Recalling the fact that $v_{1,n}\rightharpoonup 0$
in $H$ as $n\rightarrow\infty$, we have
$J'_{\lambda_n}(v_{1,n})\psi(\cdot-y_n^1)\rightarrow0$ for any $\psi\in C_0^\infty(\R^3)$, and
\begin{equation}\label{eqn:4.4}
\aligned
&{J}'_{\lambda_n}(v^1_n)\psi(\cdot-y_n^1)\\
&=(a+bA)\int_{\R^3}\nabla v_{1,n}(x+y_n^1)\nabla\psi+\int_{\R^3}V(x+y_n^1)v_{1,n}(x+y_n^1)\psi\\
&+\int_{\R^3}\lambda_nW(x+y_n^1)v_{1,n}(x+y_n^1)\psi-\int_{\R^3}f(v_{1,n}(x+y_n^1))\psi=o_n(1),
\endaligned
\end{equation}
which implies
by (\ref{eqn:deco6}) that
\begin{equation}\label{eqn:deco7}
\aligned
(a+bA)\int_{\R^3}\nabla w_1\nabla\psi&+\int_{\R^3}V(x+y_n^1)v_{1,n}(x+y_n^1)\psi\\
&+a_1\int_{\R^3}w_1\psi-\int_{\R^3}f(w_1)\psi=o_n(1).
\endaligned
\end{equation}
Set $w_\epsilon=\sqrt{|w_1|^2+\epsilon^2}-\epsilon$, $\epsilon>0$. It is clear that $w_\epsilon\rightarrow |w_1|$ in
$H$ as $\epsilon\rightarrow0$. As arguing as the previous \textbf{Claim}, we obtain (\ref{eqn:deco5}).
Let us set
\begin{equation}\label{eqn:4.11'}
v_{2,n}(x)=v_{1,n}(x)-w_1(x-y_n^1),
\end{equation}
then $v_{2,n}(\cdot+y_n^1)\rightharpoonup0$ weakly in $H$.
It follows from the Brezis-Lieb lemma that
\begin{equation}\label{eqn:4.12}
\aligned
&\|v_{2,n}\|_{s}^{s}=\|u_n\|_{s}^{s}-\|u_0\|_{s}^{s}-\|w_1\|_{s}^{s}+o(1),\,\,\text{for}\,s\in[2,6],\\
\endaligned
\end{equation}
Applying the concentration compactness principle to $\{v_{2,n}\}$, we have two possibilities:
either vanishing or non-vanishing. If vanishing occurs we have
 $$
 \sup\limits_{y\in\R^3}\int_{B_1(y)}|v_{2,n}(x)|^2\rightarrow0,
 $$
then $v_{2,n}\rightarrow0$ in
$L^s(\R^3)$ for $s\in(2,6)$, and Lemma \ref{Lem:deco} holds with $k=1$. Otherwise, $\{v_{2,n}\}$ is non-vanishing,
there exist $m'>0$ and a sequence $\{y^2_n\}\subset\R^3$ such that
\begin{equation}\label{eqn:4.04+}
\liminf\limits_{n\rightarrow\infty}\int_{B_{1}(y_n^2)}|v_{2,n}(x)|^2\geq m'>0.
\end{equation}
We repeat the arguments. By iterating
this procedure we obtain sequences of points $\{y_n^j\}\subset\R^3$ such that $|y_n^j|\rightarrow+\infty$,
$|y_n^j-y_n^i|\rightarrow+\infty$ if $i\neq j$ as $n\rightarrow+\infty$ and $v_{j,n}=v_{j-1,n}-w_{j-1}(x-y_n^{j-1})$ (like (\ref{eqn:4.11'})) with $j\geq2$
such that $v_n^j\rightharpoonup 0$ in $H$.
Based on the properties of the weak convergence, we have
$$
\aligned
(a)\quad &\|u_n\|_s^s-\|u_0\|_s^s-\sum_{i=1}^{j-1}\|w_i\|_s^s=\|u_n-u_0-\sum_{i=1}^{j-1}w_{i}(\cdot-y_n^{i})\|_s^s+o(1)\geq0,\\
(b)\quad & \text{for\,any\,} \psi\in C_0^\infty(\R^3)\, \text{with}\, \psi\geq0\,\text{and}\, i=1,...,j-1,\\
&(a+bA)\int_{\R^3}\nabla|w_i|\nabla\psi+(V_0+a_i)\int_{\R^3}|w_i|\psi
\leq \int_{\R^3}|f(w_i)|\psi.
\endaligned
$$
By the Sobolev embedding theorem and conclusion (b), we have for $i=1,...,j-1$
$$
\|w_i\|_p^2\leq S_p\int_{\R^3}((\nabla |w_i|)^2+|w_i|^2)\leq C\|w_i\|_p^p,
$$
where $S_p$ is the Sobolev constant of embedding from $H^1(\R^3)$ to $L^p(\R^3)$.
Hence, there exists $c_0>0$ independent of $w_i$ such that
$|w_i|_p^2\geq c_0$.
 Since $\{u_n\}$ is bounded sequence in $H$, conclusion (a) implies that
the iteration stop at some finite index $k$. The proof is complete.
\qed

\begin{remark}
The proof of Lemma \ref{Lem:deco} is in the spirit of the works Struwe \cite{Struwe84} and Li and Ye \cite{Li14}. It is worth of pointing out that this is the first result on decomposition of (PS) sequences
(families of approximating solutions, may be sign-changing solutions) with general energy level for Kirchhoff type equation.
We can find the decomposition of positive solution sequences with mountain pass energy level in \cite{Li14,Liu15,Tang17}, which is used to recover the compactness.

\end{remark}
Now we investigate the exponential decay property of there approximating solutions $\{u_n\}$. For notations simplicity,
in Lemma \ref{Lem:deco}, we define $y_n^0=0$, $a_0=0$ and $u_0=w_0$.
 Thus the conclusion in Lemma \ref{Lem:deco} can be restated as $|y_n^j-y_n^i|\rightarrow\infty$, $0\leq i<j\leq k$,
 $$
 \|u_n-\sum_{i=0}^{k}w_{i}(\cdot-y_n^{i})\|\rightarrow0,
 $$
 for any $\psi\in C_0^\infty(\R^3)$ with $\psi\geq0$
 \begin{equation}\label{eqn:fur-deco0}
(a+bA)\int_{\R^3}\nabla w_i\nabla\psi+(V_0+a_i)\int_{\R^3}w_i\psi
\leq \int_{\R^3}|f(w_i)|\psi,\,i=0,1,...,k.
\end{equation}

\begin{lemma}\label{Lem:zhishu}
There exists $\delta>0$ such that
\begin{equation}\label{eqn:e-estimate1}
\int_{\Omega_{R}^{(n)}}(|\nabla u_n|^2+|u_n|^2)\leq Ce^{-\delta R},
\quad \lambda_n\int_{\Omega_{R}^{(n)}}W(x)|u_n|^2\leq Ce^{-\delta R},
\end{equation}
where $\Omega_{R}^{(n)}=\R^3\setminus{\bigcup_{i=0}^kB_{R}(y_n^i)}$ and $C>0$ is independent of $n,R$.
\end{lemma}
\Proof
Using Moser's iteration to the differential inequality (\ref{eqn:fur-deco0}), we can obtain
for $i=1,...,k$
$$
\int_{\R^3\setminus B_R(0)}(|\nabla w_i|^2+|w_i|^2)\leq Ce^{-\delta R},\quad \|w_i\|_{L^\infty(\R^3\setminus B_R(0))}\leq  Ce^{-\delta R}.
$$
So by property (iii) of Lemma \ref{Lem:deco}, we have for $s\in[2,6]$
$$
\aligned
\int_{\Omega_{R}^{(n)}}|u_n|^s&\leq \|u_n-\sum_{i=0}^{k}w_{i}(\cdot-y_n^{i})\|_{L^s(\Omega_{R}^{(n)})}^s
+\sum_{i=0}^{k}\int_{\R^3\setminus B_R(0)}w_i^s\\
&\leq o_n(1)+ Ce^{-\delta R}.
\endaligned
$$
So we use Moser's iteration to prove the $L^\infty$-estimate
$$
|u_n(x)|\leq o_n(1)+ Ce^{-\delta R},\quad \text{for\,all\,}x\in \Omega_{R}^{(n)},
$$
which implies that for any $\epsilon>0$, there exist $n_0,R_0>0$ such that for $n\geq n_0$ there holds
$$
|u_n(x)|\leq \epsilon,\quad\forall x\in\Omega_{R_0}^{(n)}.
$$
Thus, in view of ($V_1$) and ($f_1$), by choosing $\epsilon, R_0$ such that for $R>R_0$, we have
\begin{equation}\label{eqn:guji--}
\int_{\R^3}\bigg(a\nabla u_n\nabla \varphi+\lambda_nW(x)u_n\varphi+\frac{V_0}{2}u_n\varphi\bigg)\leq 0,\quad\text{for\,all\,}x\in\Omega_{R}^{(n)}.
\end{equation}
For any $R>0$, define $\varphi_R$ as $\varphi_R(x)=0$ for $x\not=\Omega_{R}^{(n)}$,
$\varphi_R(x)=1$ for $x\not=\Omega_{R+1}^{(n)}$ and $|\nabla\varphi_R|\leq2$. Let
$\varphi=\varphi_R^2 u_n$, then (\ref{eqn:guji--}) can be estimated as follows:
\begin{equation}\label{eqn:guji-1}
\int_{\Omega_{R}^{(n)}}\bigg(a\nabla u_n(\varphi_R^2\nabla u_n+2u_n\varphi_R\nabla\varphi_R)+(\lambda_nW(x)+\frac{V_0}{2})u_n^2\varphi_R^2\bigg)\leq 0,
\end{equation}
which implies
\begin{equation}\label{eqn:guji-2}
\aligned
\int_{\Omega_{R}^{(n)}}(a|\nabla u_n|^2+\frac{V_0}{2}u_n^2)\varphi_R^2
&\leq C\int_{\Omega_{R}^{(n)}}|u_n\nabla u_n\varphi_R\nabla\varphi_R|\\
&\leq C\int_{\Omega_{R}^{(n)}\setminus{\Omega_{R+1}^{(n)}}}(a|\nabla u_n|^2+\frac{V_0}{2}u_n^2),
\endaligned
\end{equation}
where $C>0$ does not depend on $n, R$. From (\ref{eqn:guji-2}) we infer that
$$
\int_{\Omega_{R+1}^{(n)}}(a|\nabla u_n|^2+\frac{V_0}{2}u_n^2)
\leq \frac{C}{1+C}\int_{\Omega_{R}^{(n)}}(a|\nabla u_n|^2+\frac{V_0}{2}u_n^2).
$$
Thus, there exist $C>0$ (independent of $n, R$) and $\delta$ such that
$$
\int_{\Omega_{R}^{(n)}}(|\nabla u_n|^2+|u_n|^2)\leq Ce^{-\delta R}.
$$
Returning to (\ref{eqn:guji-1}) we also have
$$
\lambda_n\int_{\Omega_{R}^{(n)}}W(x)|u_n|^2\leq Ce^{-\delta R}.
$$
The proof is complete.
\qed

Motivated by \cite{Cerami05}, we derive a local Pohozaev-type identity which is of use in proving the convergence of solution sequences.
\begin{lemma}\label{Lem:local-p}
If $u\in E_\lambda$ solves equation ($K_{\lambda}$), then the following identity holds:
$$
\aligned
&\frac{1}{2}\int_{\R^3}t\cdot\nabla V(x)|u|^2\psi+\frac{\lambda}{2}\int_{\R^3}t\cdot \nabla W(x)|u|^2\psi\\
&=
-\frac{1}{2}\int_{\R^3}|\nabla u|^2t\cdot\nabla\psi+\int_{\R^3}t\cdot\nabla u\nabla u\cdot\nabla\psi\\
&-\frac{1}{2}\int_{\R^3}(V(x)+\lambda W(x))|u|^2t\cdot\nabla\psi+\int_{\R^3}\big(F(u)
+\frac{\lambda}{r}|u|^r\big)t\cdot\nabla\psi
\endaligned
$$
for $t\in\R^3$ and $\psi\in C_0^\infty(\R^3)$.
\end{lemma}
\Proof
Choose $\psi\in C_0^\infty(\R^3)$, $t\in\R^3$. Taking $t\cdot\nabla u\psi$
as test function in equation ($K_\lambda$) and integrating by parts,
we get the local Pohozaev-type identity. We can see \cite{Cerami05} for the details of proof.
\qed

Without loss of generality, we assume that $|y_n^1|=\min\{|y_n^i|,i=1,...,k\}$. Denote $y_n=y_n^1$ for simplicity of notations.
Borrowing from the idea in \cite{Cerami05}, we construct a sequence of cones $\mathcal{C}_n$,
having vertex $\frac{1}{2}y_n$ and generated by a ball $B_{R_n}(x_n)$ as follows:
$$
\mathcal{C}_n=\bigg\{z\in\R^3|z=\frac{1}{2}y_n+l(x-\frac{1}{2}y_n),x\in B_{R_n}(y_n),l\in[0,\infty)\bigg\},
$$
where $R_n$ satisfies
$$
\frac{\gamma}{k}\cdot\frac{|y_n|}{2}=r_n\leq R_n\leq kr_n=\gamma\cdot\frac{|y_n|}{2},\quad \gamma=\frac{1}{5(\bar{c}+1)},
$$
where $\bar{c}$ is the constant in the definition of the condition ($V_4$).
It is known in \cite{Cerami05} that the cone $\mathcal{C}_n$ has the following property:
\begin{equation}\label{eqn:po-6}
\partial\mathcal{C}_n\cap\bigcup\limits_{i=0}^{k}B_{\frac{r_n}{2}}(y_n^i)=\emptyset.
\end{equation}

\begin{lemma}\label{Lem:compact}
Let $\{u_n\}\subset E_\lambda$ be a solution sequence of ($K_\lambda$) with $\lambda=\lambda_n$. Assume that
$\|u_n\|\leq M$ for some $M>0$ independent of $n$, then, up to subsequence, there exists $u_0\in H$ such that
$u_n\rightarrow u_0$ in $H$.
\end{lemma}
\Proof
We now apply the local Pohozaev identity. Take $u=u_n$,
$t=t_n=\frac{y_n}{|y_n|}$ and $\psi=\eta\varphi_R$, where $\eta,\varphi_R\in C_0^\infty(\R^3)$ such that
$\eta(x)=0$ for $x\not\in \mathcal{C}_n$, $\eta(x)=1$ for $x\in \mathcal{C}_n$
and dist$(x,\partial\mathcal{C}_n)\geq1$, $\varphi_R(x)=1$ for $x\in B_R$, and $\varphi_R(x)=0$ for
$x\in\R^3\setminus{B_{2R}}$. By letting $R\rightarrow\infty$, we have
\begin{equation}\label{eqn:po-7}
\aligned
&\frac{1}{2}\int_{\R^3}t_n\cdot\nabla V(x)|u_n|^2\eta
+\frac{\lambda_n}{2}\int_{\R^3}t_n\cdot\nabla W(x)|u_n|^2\eta\\
=&
-\frac{1}{2}\int_{\R^3}|\nabla u_n|^2t_n\cdot\nabla\eta+\int_{\R^3}t_n\cdot\nabla u_n\nabla u_n\cdot\nabla\eta
-\frac{1}{2}\int_{\R^3}V(x)|u_n|^2t_n\cdot\nabla\eta\\&+\int_{\R^3}(F(u_n)+\frac{\lambda_n}{r}|u_n|^r)t_n\cdot\nabla\eta-
\frac{\lambda_n}{2}\int_{\R^3}W(x)|u_n|^2t_n\cdot\nabla\eta.
\endaligned
\end{equation}
From (\ref{eqn:po-6}) and the definition of $\eta$,
we see that the support of $\nabla \eta$ is contained in the domain
$\Omega=\Omega_R^{(n)}$ with $R=\frac{1}{2}r_n-1$.
In view of Lemma \ref{Lem:zhishu}, we know that the right-hand side
of (\ref{eqn:po-7}) decays exponentially, say less than $Ce^{-\delta|y_n|}$.
Observe that by Lemma 4.2 of \cite{Cerami05},
we have $t_n\cdot\nabla V\geq \frac{1}{2}\frac{\partial V}{\partial r}$ for $x\in \mathcal{C}_n$.
Besides, by the definition of $W$, we see that $\int_{\R^3}t_n\cdot\nabla W(x)|u_n|^2\eta$ is bounded uniformly for $\lambda_n$.
So the
left-hand side of can be estimated as
\begin{equation}\label{eqn:po-8}
\aligned
\frac{1}{2}&\int_{\R^3}t_n\cdot\nabla V(x)|u_n|^2\eta
+\frac{\lambda_n}{2}\int_{\R^3}t_n\cdot\nabla W(x)|u_n|^2\eta\\
&=\frac{1}{2}\int_{\R^3}t_n\cdot\nabla V(x)|u_n|^2\eta+o(1)
\geq\frac{1}{2}\inf_{x\in B_1(y_n)}\frac{\partial V(x)}{\partial r}\int_{\R^3}|u_n|^2+o(1)\\
&\geq\frac{m}{4}\inf_{x\in B_1(y_n)}\frac{\partial V(x)}{\partial r},
\endaligned
\end{equation}
where $\int_{B_1(y_n)}u_n^2dx\geq m>0$.  Thus,
 together (\ref{eqn:po-7}) and (\ref{eqn:po-8}), we obtain
$$
\frac{m}{4}\inf_{x\in B_1(y_n)}\frac{\partial V(x)}{\partial r}\leq Ce^{-\delta|y_n|},
$$
which contradicts with (V$_3$). Thus $k=0$ and by Lemma \ref{Lem:deco} (iii), we have
$u_n\rightarrow u_0$ in $H$.
\qed

In view of Lemma \ref{Lem:compact}, $u_{0}$ is a nontrivial solution of problem (K).
Actually we have proved the following fact.
\begin{proposition}\label{prop:convergence}
Assume $\{u_{\lambda}\}_{\lambda\in(0,1]}$ satisfies $I'_{\lambda}(u_{\lambda})=0$ and
$c_{\lambda}=I_{\lambda}(u_{\lambda}) \in [m_{1},m_{2}]$, then
there exists $u_{0}\in H\setminus\{0\}$ such that on a sequence
$\{\lambda_{n}\}$ tending to zero, it holds
$$u_{\lambda_{n}} \to u_{0} \quad \text{in } H,
\quad c_{\lambda_{n}}\to c_{0}, \quad I(u_{0}) = c_{0} \quad \text{ and } \quad I'(u_{0}) = 0.$$
\end{proposition}

Based on Proposition \ref{prop:convergence}, we are now able to give the

\subsection{Proof of Theorem \ref{Thm:existence}}

Define the set of solutions
$$
\mathcal{S}:=\{u\in H\setminus\{0\}:\,I'(u)=0\}
$$
that, for what we have proved, is nonempty.
For  $u\in \mathcal{S}$, by Sobolev's inequality, for any $\varepsilon>0$ there exists $C_\varepsilon>0$ such that
$$
\aligned
 \|u\|^2+b\|\nabla u\|_2^4\leq \varepsilon\int_{\R^3}u^2+C_\varepsilon\int_{\R^3}|u|^6
\endaligned
$$
which implies  that $\mathcal S$ is bounded away from zero. Besides, we can also see from the above inequality that
$\|\nabla u\|_2^2\geq C$ for all $u\in \mathcal{S}$.
By recalling (\ref{eqn:fun1})-(\ref{eqn:fun3}),
 there exists some $C>0$ satisfying $I(u)\geq C\|\nabla u\|_2^2$ for all $u\in \mathcal{S}$.
So we infer that
$$
c_*:=\inf\limits_{u\in\mathcal{S}}I(u)>0.
$$
Choose finally a minimising sequence
$\{u_n\}\subset \mathcal{S}$ so that $I(u_n)\rightarrow c_*$.
Similarly to Lemma \ref{Lem:local} we know that
 $\{u_n\}$ is bounded in $H$. Like the modified functional $I_\lambda$, we can also prove
 some facts for solution sequence $\{u_n\}$ of $I$
 corresponding to Lemmas \ref{Lem:deco}-\ref{Lem:compact}.
As a consequence,
there exists $u_{*}\in H$ so that $u_n\rightarrow u_{*}$ in $H$ and $I'(u_{*})=0$.
Then $u_{*}$ is a ground state solution of (K).

\s{Multiplicity}\label{sec:critical}
In this section, we are attempt to use the perturbation approach together with the Symmetric Mountain-Pass
theorem to prove the existence of infinitely many high energy solutions to problem (K).
\subsection{Proof of Theorem \ref{Thm:duochjie}}
We recall that
$I_{\lambda}$ belongs to $C^{1}(E,\R)$. Denote $\mathcal{B}_{R}$ by the ball of radius $R>0$ of $E$.
Choose a sequence of finite dimensional subspaces $E_j$ of $E$ such that $\dim E_j=j$
and $E_j^{\bot}$ denotes the orthogonal complement of $E_j$.
We define $\partial \mathcal{P}$ by
$$
\aligned
\partial \mathcal{P}:=\bigg\{u\in E\setminus\{0\}\,\bigg|\,&\frac{(\mu+2)a}{2\mu}\int_{\R^3}|\nabla u|^2
+\frac{2+3\mu}{2\mu}\int_{\R^3}V(x)u^2+\frac{1}{2}\int_{\R^3}(\nabla V(x),x)u^2\\
&+\frac{(2+3\mu)\lambda}{2\mu}\int_{\R^3}W(x)u^2+\frac{\lambda}{2}\int_{\R^3}(\nabla W(x),x)u^2
+\frac{(\mu+2)b}{2\mu}\bigg(\int_{\R^3}|\nabla u|^2\bigg)^2\\
&=\int_{\R^3}(\frac{1}{\mu}f(u)u+3F(u))+\frac{(r+3\mu)\lambda}{\mu r}\int_{\R^3}|u|^r
\bigg\}.
\endaligned
$$
Recalling assumption $(V_4)$ and (\ref{eqn:differe}), it follows from Sobolev's inequality that
for any $\varepsilon>0$, there exists $C_\varepsilon>0$ such that
\begin{equation}\label{eqn:many1}
\aligned
&\frac{(\mu+2)a}{2\mu}\int_{\R^3}|\nabla u|^2+
\frac{2+3\mu}{2\mu}\int_{\R^3}V(x)u^2\\
&\leq \int_{\R^3}(\frac{1}{\mu}f(u)u+3F(u))+\frac{(r+3\mu)}{\mu r}\int_{\R^3}|u|^r\\
&\leq \varepsilon\int_{\R^3}|u|^2+C_\varepsilon\int_{\R^3}|u|^6,\quad\forall u\in \partial \mathcal{P}\cap E_j^{\bot},
\endaligned
\end{equation}
which implies that there exists $m_3>0$ independent of $\lambda$ such that $\|\nabla u\|_2^2\geq m_3$.
For any $u\in \partial \mathcal{P}\cap E_j^{\bot}$, using the definition of $I_\lambda$ and (\ref{eqn:many1}), we arrive at
\begin{equation}\label{eqn:many2}
\aligned
I_\lambda(u)&\geq a\frac{3\mu-2}{8\mu}\int_{\R^3}|\nabla u_n|^2+\frac{\mu-2}{8\mu}b\left(\int_{\R^3}|\nabla u_n|^2\right)^2
+\lambda\frac{r-\mu}{4\mu r}\int_{\R^3}|u_n|^r\\
&\geq a\frac{3\mu-2}{8\mu}m_3+\frac{\mu-2}{8\mu}bm_3^2=:\delta.
\endaligned
\end{equation}
Moreover, we can choose $R_j>0$ such that $I_\lambda(u)<0$ for $u\in E_j\cap\partial \mathcal{B}_{R_j}$.
Actually,  such an $R_{j}$ can be found by the fact that
in the proof of (2) of Lemma \ref{Lem:MP1} the element $e\in C_0^\infty(\R^3)$ is arbitrary.
Note that $R_{j}$ does not depends on $\lambda$, that is to say,
$$
\forall \lambda\in(0,1] : I_{\lambda} (u)<0 \quad\text{ for any } u \in E_j\cap\partial \mathcal{B}_{R_j}.
$$
Thus, the functional $I_{\lambda}$
satisfies all the assumptions of the Symmetric Mountain Pass Theorem,
and we define the minimax values
$$
 c_\lambda(j)=\inf\limits_{B\in\Gamma_j}\sup\limits_{u\in B} I_\lambda(u)
 $$
 where
$$
\Gamma_j=\bigg\{B=\phi(E_j\cap \mathcal{B}_{R_j}) | \phi\in C(E_j\cap \mathcal{B}_{R_j},E),
\,\phi \,\,\text{is odd},\ \phi=\textrm{Id} \ \text{on} \ E_j\cap \partial \mathcal{B}_{R_j}\bigg\}.
$$
It is easy to prove that the following intersection property holds (see \cite[Proposition 9.23]{Rabinowitz86}):
for $B\in\Gamma_j$,
$$
B\cap\partial \mathcal{P}\cap E_j^{\bot}\not=\emptyset,
$$
which implies by (\ref{eqn:many2}) that $c_\lambda(j)>\delta>0$.
 For any fixed $j$, by the definition of $ c_\lambda(j)$, we have, in view of (2) of Lemma \ref{Lem:MP1},
$$
\aligned
 c_\lambda(j)&\leq\sup\limits_{u\in E_j\cap \mathcal{B}_{R_j}} I_\lambda(u)\\
 &\leq \sup\limits_{u\in E_j\cap \mathcal{B}_{R_j}} \bigg
 \{ C_{1} \|u\|_{E}^2+ C_{2}\| u\|_{E}^{4}\bigg\}:=C_{R_j},
\endaligned
$$
where $C_{R_j}$ is indeed independent of $\lambda\in(0,1]$ and $\|\cdot\|_{E}$ is any norm
in $E_{j}$.
Based on the above arguments, one has $c_\lambda(j)\in[\delta, C_{R_j}]$. Using again Lemmas \ref{Lem:deco}-\ref{Lem:compact} and
Proposition \ref{prop:convergence}, we know that there exists $u_{0}(j)\in H\setminus\{0\}$
such that on a sequence $\lambda_{n}\to 0^{+}$,
$$u_{\lambda_{n}}(j) \to u_{0}(j) \quad \text{in } H,
\quad c_{\lambda_{n}}(j)\to c_{0}(j)\geq\delta, \quad I(u_{0}(j)) = c_{0}(j) \quad \text{ and } \quad I'(u_{0}(j)) = 0,$$
that is, $u_{0}(j)$ is a nontrivial solution of problem (K).

Once we show that $c_0(j)\rightarrow+\infty$ as $j\to+\infty$,
problem (K) has infinitely many bounded state solutions and the proof
of Theorem \ref{Thm:duochjie} is finish.

Now we give an estimate on $I_\lambda$ as follows
$$
\aligned
I_\lambda(u)=&I(u)+\frac{\lambda}{2}\int_{\R^3}W(x)u^2-\frac{\lambda}{r}\int_{\R^3}|u|^r\\
\geq&\frac{1}{2}\int_{\R^3}(|\nabla u|^2+V(x)u^2)-\frac{1}{r}\int_{\R^3}|u|^{r}:=L(u).
\endaligned
$$
Define the  set $\partial\Theta\subset H$ by
$$
\partial\Theta:=\bigg\{u\in H\setminus\{0\}: \int_{\R^3}(|\nabla u|^2+V(x)u^2)=\int_{\R^3}|u|^{r}\bigg\},
$$
which is the Nehari manifold associated to energy functional $L$,
which, by classical arguments, is bounded away from zero and homeomorphic to the unit sphere.
Then, for $B\in \Gamma_j$, an easy modification of the proof of \cite[Proposition 9.23]{Rabinowitz86}
shows that an intersection property holds so that $\gamma(B\cap\partial\Theta)\geq j$, for all $j\in \mathbb N$.
%
Here $\gamma(\cdot)$
denotes the Krasnoselski genus of a symmetric set. Hence,
$$
 c_{\lambda}(j)=\inf\limits_{B\in\Gamma_j}\sup\limits_{u\in B}  I_\lambda(u)
 \geq\inf\limits_{A\subset\partial\Theta,\gamma (A)\geq j}\sup\limits_{u\in A} L(u)
:=b(j).
$$
It is not hard to verify that the functional $J$ is bounded below on $\partial\Theta$.
Moreover, We observe that
the boundedness of the Palais-Smale sequence is
easy to verify for functional $J$. As a result,
with some suitable modification, the arguments of functional $I_{\lambda}$ are still
valid for $J$ without any perturbation. So,
$J$ satisfies the Palais-Smale condition. Then the Ljusternick-Schnirelmann theory
guarantees that  $b(j)$ are diverging critical values for $J$.
 Therefore,
$$
 c_0(j)=\lim\limits_{\lambda\rightarrow0^+}c_{\lambda}(j)\geq b(j)\rightarrow+\infty, \ \ \text{as } j\to+\infty.
$$
That is to say, problem \eqref{K} has infinitely many higher energy solutions. The proof is complete.
\qed

\end{document}